\documentclass[11pt]{AAS}

\usepackage{bm}
\usepackage{amsmath}
\usepackage{subfigure}
\usepackage{overcite}
\usepackage{footnpag}			      

\usepackage{todonotes}
\usepackage{rotating}

\PaperNumber{23-130}

\begin{document}

\title{Tracing Position in the Regime of the Restricted Three-Body Problem to a Halo Orbit}

\author{Hailee E. Hettrick\thanks{Ph.D. Candidate \& Draper Scholar, Aeronautics \& Astronautics, Massachusetts Institute of Technology, The Charles Stark Draper Laboratory, 70 Vassar St, Cambridge, MA 02139.},  
Begum Cannataro\thanks{Principal Member of Technical Staff, Autonomy \& Real-Time Planning, The Charles Stark Draper Laboratory, 555 Technology Square, Cambridge, MA 02139.},
\ and David W. Miller\thanks{Professor Post-Tenure, Aeronautics \& Astronautics, Massachusetts Institute of Technology, 70 Vassar St, Cambridge, MA 02139.}
}
\maketitle{}

\begin{abstract}
Driven by the desire to find positions that satisfy keepout constraints for a space-based telescope mission, this work develops a process for tracing a point in space in the regime of the restricted three-body problem to a halo orbit, characterized by its out-of-plane amplitude, and its position on that halo orbit, denoted by the halo orbit time. This process utilizes third-order solutions from the Lindstedt-Poincar\'e method, which have been partially inverted to expect a point in space as an input. Three different methodologies that use these partially inverted expressions are presented. Results are produced for 1,000 randomly selected points using all three methods and are compared to truth. Ultimately, the method that employed two distinct accuracy metrics yielded the most accurate results for the dataset.
\end{abstract}

\section{Introduction}
In certain solar system dynamical regimes, nonlinear dynamics exist representing the balance of gravitational and centrifugal forces on a negligibly small body imparted by two massive celestial bodies. The motion of the three bodies is referred to as the restricted three-body problem, and a prefix is attached which describes the assumed motion of the massive bodies. Of particular interest is the circular restricted three-body problem (CR3BP). CR3BP has five equilibrium points (also known as libration or Lagrange points); these points consist of three collinear equilibrium points\cite{euler1767motu} and two off-axis equilibrium points\cite{lagrange1772essai}. These equilibrium points are typically denoted as $L_\#$ corresponding to their position.

Periodic orbits exist around these Lagrange points; of particular interest is the halo orbit. Farquhar discovered a three-dimensional periodic solution to the Earth-Moon CR3BP that he called a ``halo orbit" because of how the orbit looked like a halo around the Moon as viewed from Earth \cite{farquhar1968control}. This three-dimensional periodic solution is a trajectory corresponding to high-precision initial conditions, which are found numerically via a differential corrector \cite{breakwell1979halo}. For this scheme to quickly converge to initial conditions corresponding to a periodic orbit, the initial guess must be in the neighborhood of the solution. Richardson provided a method for solving a third-order approximation to these initial conditions given a desired amplitude of the orbit; the resulting third-order approximation serves as an initial guess for a differential corrector \cite{richardson1980analytic} \cite{richardson1980halo}. Howell generalized the work of her predecessors, which was focused on the Earth-Moon system, to any system characterized by the mass ratio of the two primary bodies in the CR3BP \cite{howell1984three}.

Research into halo orbits then led to the investigation of how to use invariant manifolds, which serve as natural pathways to periodic orbits about a Lagrange point. With these natural pathways, a spacecraft can move from one periodic orbit to another along such a manifold without expending fuel. These manifolds were not well understood in space mission planning until the end of the twentieth century. G\'omez et al. identified a halo orbit about Sun-Earth $L_1$ fitting desired mission parameters and then globalized the corresponding stable manifold terminating at the halo orbit via backward integration \cite{gomez1991study}. G\'omez et al. further considered these natural pathways from Earth to a halo orbit about Sun-Earth $L_1$ by accounting for the Moon's influence on the stable manifold (i.e., considering the system Sun-Earth+Moon rather than only the Sun-Earth system) \cite{gomez1991moon} \cite{gomez1993study}. In this manner, Gomez et al. introduced the construction of transfer trajectories to a periodic orbit about a Lagrange point using invariant manifold theory. Howell, Mains, and Barden expanded on this strategy by adding a differential corrector in tandem to generate transfers from a given Earth parking orbit \cite{howell1994transfer}. Furthermore, Barden developed a modified differential corrector scheme to generate transfers from Earth-to-halo orbit and halo-to-halo \cite{barden1994using}. 

The first spacecraft mission to utilize the periodic orbits (specifically, a halo orbit) about a Lagrange point (Sun-Earth $L_1$ and later Sun-Earth $L_2$) was the third International Sun-Earth Explorer (ISEE-3), launched in 1978 \cite{farquhar2001flight}. The successful use of the halo orbit in ISEE-3's mission led to the use of periodic orbits (halo or Lissajous orbits) for several more missions: SOHO (1996) \cite{domingo1995soho}, ACE (1997) \cite{stone1998advanced}, MAP (2001) \cite{bennett2003microwave}, Genesis (2001) \cite{lo2001genesis}, Planck (2009) \cite{passvogel2010planck}, Herschel Space Observatory (2009) \cite{pilbratt2010herschel}, Chang'e 2 (2010) \cite{wu2012pre}, THEMIS-B and THEMIS-C (2011) \cite{broschart2009preliminary}, Gaia (2013) \cite{prusti2016gaia}, DSCOVR (2015) \cite{roberts2015early}, LISA Pathfinder (2015) \cite{mcnamara2008lisa}, Chang'e 4 (2018) \cite{jia2018scientific}, Chang'e 4 Relay (2018) \cite{gao2019trajectory}, Specktr-RG (2019) \cite{kovalenko2019orbit}, and JWST (2021) \cite{gardner2006james}.

Lunar Gateway and the Habitable Exoplanet Observatory (HabEx) will build on the success of these past missions that exploit halo orbits. Lunar Gateway is a prominent future mission intending to make use of a periodic orbit about a Lagrange point (specifically, a near-rectilinear halo orbit about Earth-Moon $L_2$) and is a cornerstone of NASA's Artemis program \cite{mars_2016}. A 12U CubeSat, Cislunar Autonomous Positioning System Technology Operations and Navigation Experiment (CAPSTONE), will serve as an orbital precursor for Lunar Gateway to verify the orbital stability of Lunar Gateway's future intended orbit \cite{hall_2020}. HabEx was proposed in the 2020 Astronomy and Astrophysics Decadal Survey and aims to use a space-based telescope and an external occulter in formation flight about Sun-Earth $L_2$ to observe exoplanets. This type of mission requires the external occulter to be within a defined tolerance of axial separation from the telescope during imaging to properly observe the target exoplanet.

HabEx's mission stresses the importance of maintaining the formation flight of the two spacecraft during imaging. Since both spacecraft will be placed in the regime of the restricted three-body problem, utilizing halo orbits would increase science yield and minimize fuel expenditures. Yet, successful observations of target exoplanets face constraints posed by bright celestial bodies, primarily the Sun and Earth. This approach facilitates the determination of constraint-satisfying positions in space at specific moments in time. However, the evaluation of the orbits associated with those constraint-satisfying positions is still pending.  Therefore, this work develops a process that determines the halo orbit and the location along the halo orbit that coincides with a point in space satisfying these conditions.

\section{Problem Formulation}
Tracing a point in space in the regime of the restricted three-body problem to a halo orbit first requires a discussion of CR3BP and techniques that yield approximate analytical closed-form solutions. CR3BP describes the motion of an object with negligible mass under the influence of the competing gravitational pull of two massive celestial bodies (denoted as \emph{primaries}), which are assumed to orbit each other circularly \cite{koon2000dynamical}. Traditionally, the problem is considered in a rotating frame, such that the primaries appear stationary along the $x-$axis and the $x-y$ plane is the ecliptic. Additionally, the distance between the primaries is normalized to be unity, and the angular velocity of the rotating frame is also normalized to unity. The masses of the primaries, $m_1$ and $m_2$, are used to form the mass ratio $\mu = m_2/(m_1+m_2)$. Considering the larger primary to be located at $(-\mu,0,0)$ and the smaller primary to be at $(1-\mu,0,0)$, the kinetic and potential energies of the massless body are found.

\begin{equation}
\begin{aligned}
	K &= \frac{1}{2}\biggr( (\dot{x} - y)^2 + (\dot{y}+x)^2 \biggr) \\
	U &= - \frac{1-\mu}{\sqrt{(x+\mu)^2+y^2+z^2}} - \frac{\mu}{\sqrt{(x-1+\mu)^2+y^2+z^2}}
\end{aligned}
\label{eqn:energies}
\end{equation}

\begin{figure}[htb]
	\centering\includegraphics[width=3.5in]{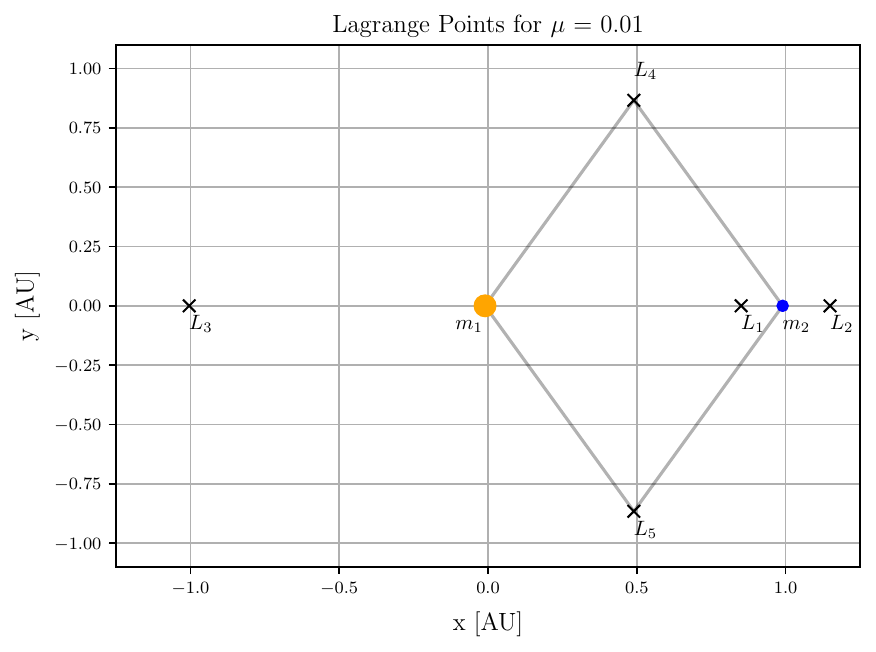}
	\caption{Equilibrium points of the restricted three-body problem}
	\label{fig:Lpts}
\end{figure}

Defining the Lagrangian as $L = K-U$ and applying the Euler-Lagrange equations yields the massless body's equation of motions as
\begin{equation}
	\begin{aligned}
		\ddot{x} - 2\dot{y} &= \frac{\partial \bar{U}}{\partial x} \\
		\ddot{y} + 2\dot{x} &= \frac{\partial \bar{U}}{\partial y} \\
		\ddot{z} &= \frac{\partial \bar{U}}{\partial z}
	\end{aligned}
\label{eqn:eom}
\end{equation}
where $\bar{U} = \frac{1}{2}(x^2+y^2) - U $. Using Eqs.~\eqref{eqn:energies} and ~\eqref{eqn:eom}, the five Lagrange points for the system can be determined via Newton's method and are dependent on the mass ratio ($\mu$) of the system. Figure \ref{fig:Lpts} illustrates the five Lagrange points for a massless body under the influence of the primaries with $\mu = 0.01$.

As mentioned previously, there exist periodic orbits about these Lagrange points that represent the balance of the gravitational forces with the centrifugal forces. Linearizing around one of the collinear equilibrium points and evaluating the eigendata indicates that the equilibrium point has center $\times$ center $\times$ saddle stability \cite{sanchez2020towards}. This corresponds to the equilibrium point having a one-dimensional stable manifold, a one-dimensional unstable manifold, and a four-dimensional center manifold \cite{barden1998fundamental}. Furthermore, the center manifold contains periodic and quasi-periodic motion. Barden and Howell identified the following types of periodic and quasi-periodic motion in the center manifold \cite{barden1998fundamental}:
\begin{enumerate}
	\item Lyapunov orbit: periodic solution in the plane of motion of the primaries
	\item Nearly vertical orbit: periodic solution dominated by out-of-plane motion
	\item Lissajous trajectory: three-dimensional quasi-periodic motion containing frequencies from both Lyapunov and nearly vertical orbits and laying on a two-dimensional torus
	\item Halo orbit: three-dimensional periodic motion resulting from a bifurcation of the Lyapunov orbit
	\item Quasi-periodic tori enveloping the halo orbits
\end{enumerate}

The work described herein focuses on the halo orbit. Typically, solutions for halo orbits are found by first evaluating a good initial condition guess using Richardson's third-order approximate solution for periodic motion \cite{richardson1980halo} \cite{richardson1980analytic}. Richardson's method uses successive approximations in conjunction with the Lindstedt-Poincar\'e method, which is detailed by Koon et al. \cite{koon2000dynamical}. Effectively, the user selects the desired out-of-plane amplitude of the orbit, and Richardson's method produces an initial condition corresponding to that amplitude.

Given an approximation of the initial conditions for a halo orbit using Richardson's method, the differential corrector modifies the initial state to guarantee that the periodicity of the orbit is achieved within some defined tolerance by using a variational relationship and iteration. This differential correction scheme is commonly used, and the algorithm is outlined particularly well by Howell \cite{howell1984three}, Soto \cite{soto2020orbital}, and Sanchez \cite{sanchez2020towards}. The refined initial condition yielded by the differential corrector is used to solve the differential equations for a desired length of time, producing a halo orbit.

 \begin{figure}[htb]
	\centering\includegraphics[width=4.5in]{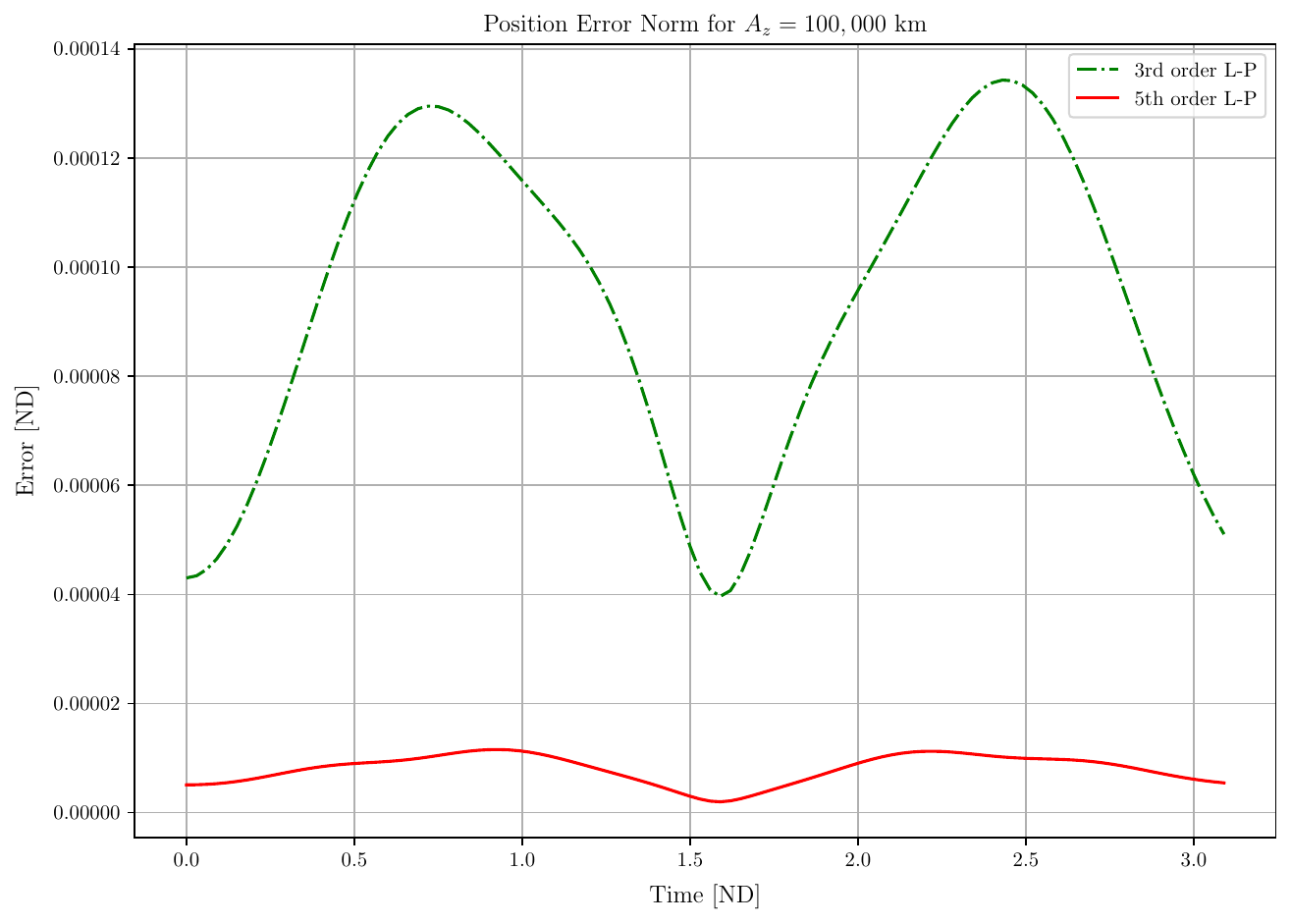}
	\caption{Position error norms for third- and fifth-order solutions }
	\label{fig:lpSol}
\end{figure}

 \begin{figure}[htb!]
	\centering\includegraphics[width=4.5in]{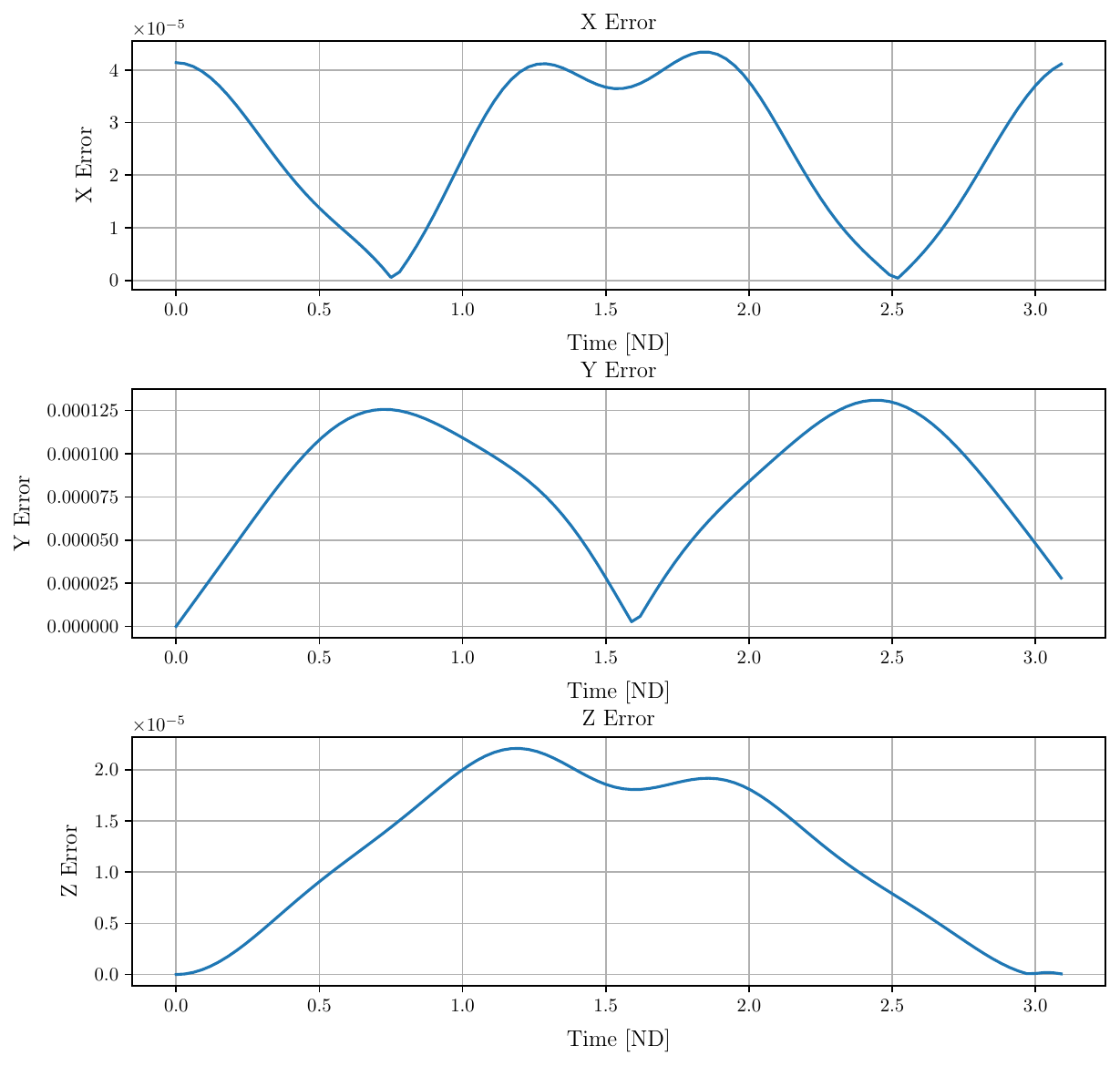}
	\caption{Position errors for third-order solution}
	\label{fig:threeCoordErrors}
\end{figure}

Rather than using the differential equations, this work uses the iterative technique of applying the Lindstedt-Poincar\'e method, detailed by Jorba and Masdemont \cite{Jorba1998}, to produce closed-form approximate analytical solutions. Ultimately, closed-form solutions are of the form 
\begin{equation}
\begin{aligned}
	x(t,A_z) &= \sum_{i,j=1}^n \biggr(\sum_{|k| \leq i+j} x_{ijk} \cos (k\theta) \biggr) A_x(A_z)^i A_z^j \\
	y(t,A_z) &= \sum_{i,j=1}^n \biggr(\sum_{|k| \leq i+j} y_{ijk} \sin (k\theta) \biggr) A_x(A_z)^i A_z^j \\
	z(t,A_z) &= \sum_{i,j=1}^n \biggr(\sum_{|k| \leq i+j} z_{ijk} \cos (k\theta) \biggr) A_x(A_z)^i A_z^j 
\end{aligned}	
\label{eqn:lp}
\end{equation}
where $n$ is the order of the solution desired, $A_z$ is the out-of-plane amplitude, and $A_x$ is the in-plane amplitude. $A_x$ is a function of $A_z$ since halo orbits have equivalent in-plane and out-of-plane frequencies. To illustrate the nature of the position error norm time history, the errors associated with third- and fifth-order solutions of a halo orbit with $A_z = 100,000$ km are shown in Figure \ref{fig:lpSol}. Additionally, errors for each position coordinate for $A_z = 100,000$ km are shown in Figure \ref{fig:threeCoordErrors}. Note that although the error seems bounded and oscillatory, it is growing with time due to the hyperbolic nature of halo orbits.

\section{Inverse Functions of Position} 
\begin{figure}[htb]
	\centering\includegraphics[width=5.5in]{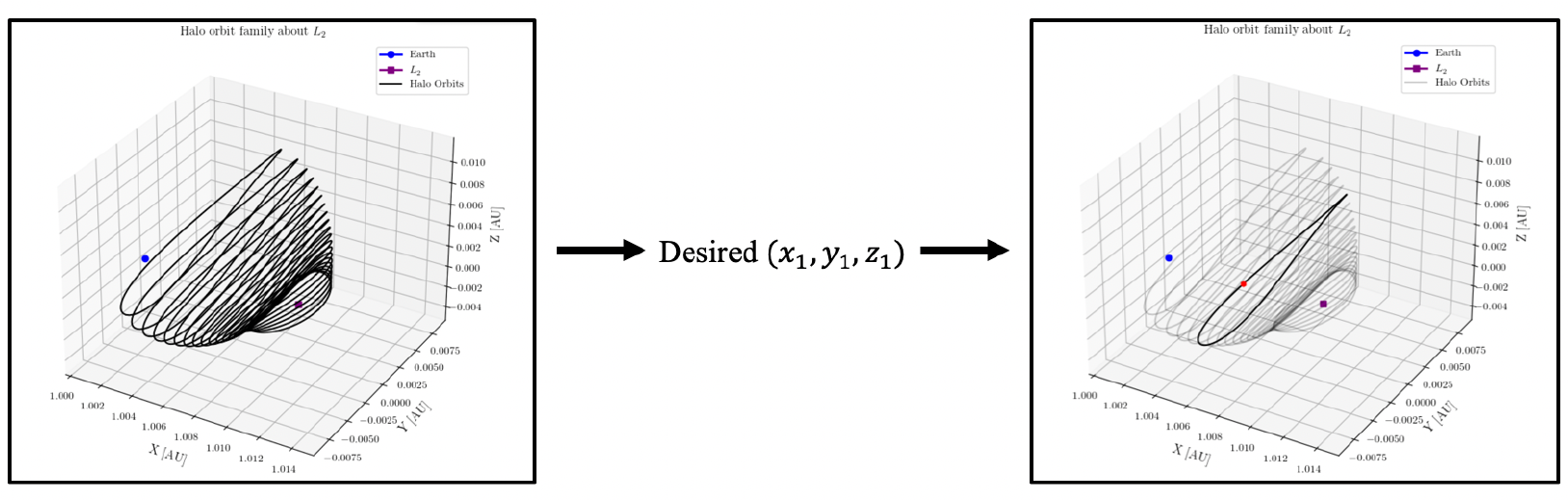}
	\caption{Find the halo orbit associated with the desired point (red) from a family of halo orbits}
	\label{fig:process}
\end{figure}

The goal of this work is to re-formulate the Lindstedt-Poincar\'e solutions in Eq.~\eqref{eqn:lp} in terms of some known point in space, denoted as $(x_1,y_1,z_1)$, to determine the halo orbit the point is on. A point on a halo orbit is characterized by the halo orbit's out-of-plane amplitude $A_z$ and time in the halo orbit's period $t$. Figure \ref{fig:process} illustrates that given a dense family of halo orbits around Lagrange points, some desired point $(x_1,y_1,z_1)$ in the regime of CR3BP can be identified as a point on a halo orbit.

Using Eq.~\eqref{eqn:lp} with $n=3$, the $x$ and $z$ equations are selected for finding inverse functions of $t$ and $A_z$ since they are both functions of cosine. First, an expression for $t$ is found using the $z$ equation and assuming some known value of $z=z_1$. This expression is found via Mathematica by solving
\begin{equation}
	z_1 = \biggr(z_{110} A_x A_z + z_{011} A_z \cos{(\omega t)} + A_x  A_z z_{112} \cos{(2 \omega t)} + z_{213} A_x^2 A_z \cos{(3 \omega t)}\biggr) \biggr| 1 - ||L_{\#}|| \biggr|
\end{equation}
for $t$. Recall that $A_x$ is defined as a function of $A_z$ -- the dependency is dropped for readability $A_x$ is used here for simplicity. Additionally, note that in these expressions, it is assumed that $A_x$ and $A_z$ have been scaled to be non-dimensional and that $L_\#$ indicates the Lagrange point of interest. The resulting expression yields $t$ as a function of $z_1$ and $A_z$ which is too lengthy to reproduce here. Therefore, the expression is referred to as $t = f(z_1, A_z)$. Plugging this expression into the $x$ equation in Eq.~\eqref{eqn:lp} produces
\begin{equation}
\begin{aligned}
x(z_1,A_z) = \biggr( ||L_{\#}|| -\mu \biggr) &+ \biggr(x_{020}A_z^2 + x_{200}Ax^2 + A_x \cos{(\omega f(z_1,A_z))} + \\
&(x_{022}A_z^2 + x_{202} A_x^2) \cos(2\omega f(z_1,A_z)) +\\
& (x_{123} A_x A_z^2 + x_{303} A_x^2 ) \cos(3 \omega f(z_1,A_z))  \biggr) \biggr| 1 - ||L_{\#}|| \biggr|.
\end{aligned}
\end{equation}

The dependency on $t$ has thus been removed. Assuming a known point $x_1$, the above would ideally be solved for $A_z$ to produce an expression $A_z = h(x_1,z_1)$. Unfortunately, Mathematica requires more than 32GB of memory to solve the expression $x_1 = g(z_1, A_z)$ for $A_z$. A solution for the expression $A_z = h(x_1,z_1)$ would be substituted into the $t$ expression so that expressions for both $t$ and $A_z$ only depend on $x_1$ and $z_1$. However, since a solution for the aforementioned expression is not available, a method is developed to minimize $x_1 - g(z_1, A_z)$ in terms of $A_z$ to then solve for $t$. 

 Due to the construction of the $t$ expression as a function of $\arccos$, the inverse function solves for values $0 \leq t \leq \frac{T}{2}$ where $T$ indicates the period of the halo's orbit. Since the halo orbit is unknown at this juncture, an estimate for the period occurs when the sign of $y_1$ indicates the point is in the second half of the orbit's period. For both Northern and Southern halo families, negative $y_1$ indicates that $t$ should be in the region of $[0, \frac{T}{2}]$; positive $y_1$ indicates that $t$ should be in the region of $(\frac{T}{2},T]$.

Since halo orbits' periods are usually around $\pi$ (when non-dimensionalized) and $y = 0$ at $t=\frac{T}{2}$, an initial guess of $\frac{T}{2} = \frac{\pi}{2}$ is used to initialize the process of finding when the Lindstedt-Poincar\'e solution produces $y = 0$. The $\frac{T}{2}$ guess is varied in the neighborhood of $\frac{\pi}{2}$ until $y=0$. While this approximation of the halo orbit's period is necessary to find $t$ in the correct half of the orbit, it does introduce another source of error.

\section{Employing the Inverse Functions}
Three different methods are developed and presented which, given some point $(x_1,y_1,z_1)$, produce the corresponding halo characteristics $(t, A_z)$. Method 1 prioritizes minimizing $x_1 - g(z_1, A_z)$, agnostic to the accuracy of the results for $y_1$ and $z_1$. Methods 2 and 3 shift the focus to minimizing the position error norm. Originally, Method 1 was the sole methodology used since minimizing only $x_1 - g(z_1, A_z)$ does not require the evaluation of the Lindstedt-Poincar\'e solutions (which is necessary for calculating the position error norm). However, Method 1's shortcomings, which will be discussed below, indicated that the position error norm must be considered. Each method is discussed in turn and results are shown and compared.

A sample of 1,000 randomly selected halo orbits is found to measure how well the different methods find the halo orbit on which a point in space belongs. From the sample of random orbits, a random index from each trajectory indicates a corresponding position $(x_1,y_1,z_1)$. In this manner, a random selection of positions in the regime of CR3BP is inputted into the process described above to yield the corresponding $(t, A_z)$ values. After computing the corresponding $(t, A_z)$ pair, they are used as inputs into the Lindstedt-Poincar\`e solutions, Eq. ~\eqref{eqn:lp}, to see how closely $(x_2,y_2,z_2)$ matches the original point. The results described herein used Northern halo families about Sun-Earth $L_2$.

\subsection{Finding $(t,A_z)$: Method 1}
Having found the functions $x = g(z_1,A_z)$ and $t = f(z_1,A_z)$, the primary process of finding $(t,A_z)$ corresponding to $(x_1,y_1,z_1)$ is executed as follows.
\begin{enumerate}
\item Initialize a range of $A_z$ values between $100$ km and $1.0 \times 10^6$ km.	
\item Initialize a variable $x_0$ to a large magnitude.
\item In a for loop,
	\begin{enumerate}
	\item Plug in $A_z[i]$ in the expression $x_1 - g(z_1,A_z)$. If the resulting magnitude is less than $x_0$, the index $i$ is stored and the value of $x_0$ is updated.
	\end{enumerate}
	\item Upon completing the for loop, the value of $x_0$ is compared against a user-defined tolerance equal to $10^{-7}$ here. 
	\begin{enumerate}
	\item If $x_0 \leq tol$, the process has found an acceptable value of $A_z$ corresponding with the stored index. The method then proceeds to find the associated value of $t$.
	\item If  $x_0 > tol$, the value of $A_z$ associated with the smallest $x_0$ is used to define a new range of $A_z$ values centered around it. This range gets smaller upon each iteration. If the user-defined number of iterations is exceeded, the process exits and indicates that no $A_z$ value could be found.
	\end{enumerate}
	\item The value of $A_z$ that satisfied the tolerance requirement is plugged into $t = f(z_1,A_z)$. If the sign of $y_1$ is positive, then $t$ is adjusted to be in the second half of orbit as discussed above.
\end{enumerate}
The outputs of Method 1 are either values of $(t, A_z)$ or empty arrays indicating that the method could not find a solution that satisfied the tolerance requirements.

\subsubsection{Results for Method 1}
\begin{figure}[h!]
	\centering\includegraphics[width=0.8\textwidth]{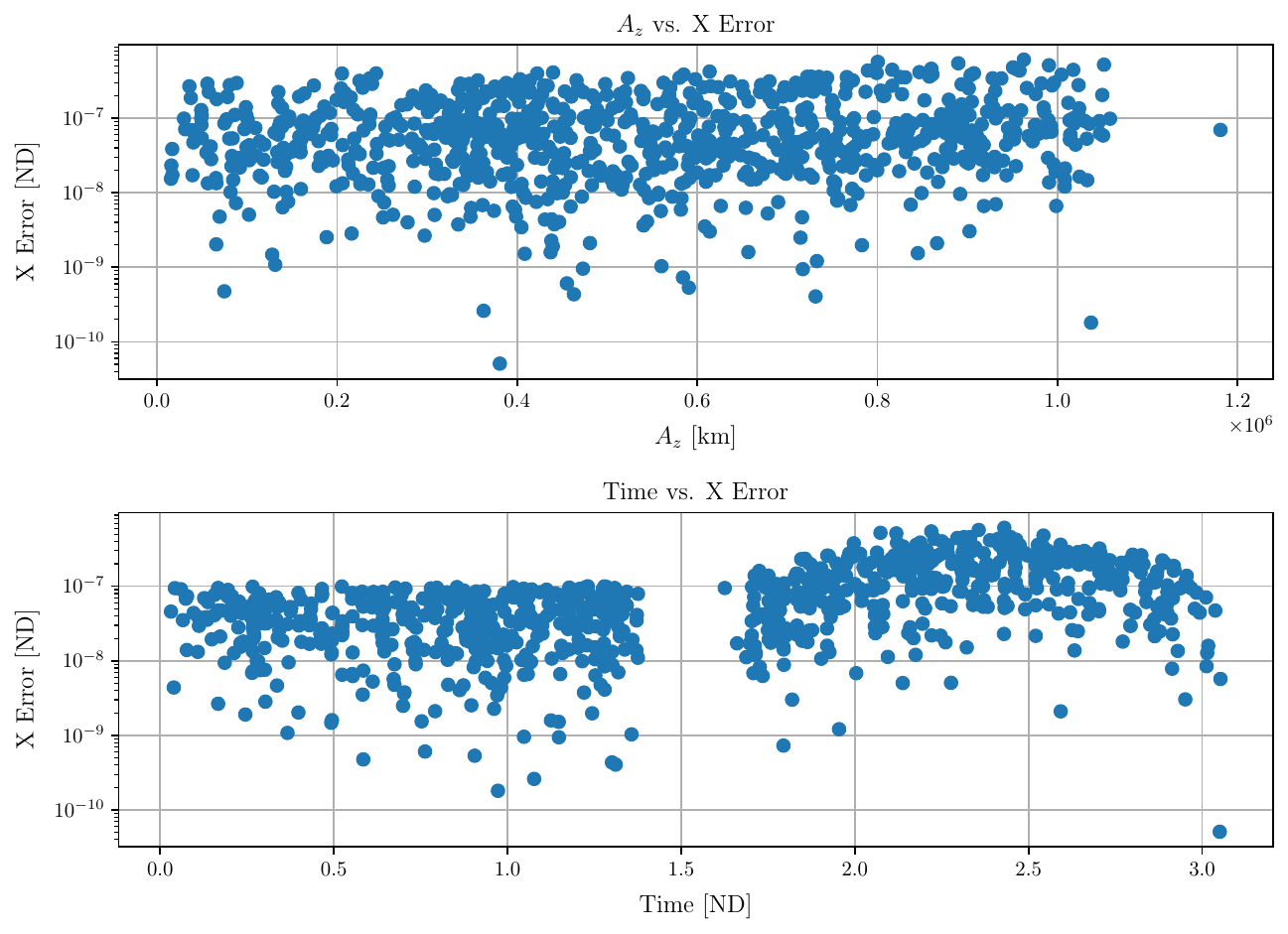}
	\caption{$x$ error compared to time and $A_z$ from Method 1}
	\label{fig:xNorm}
\end{figure}
To evaluate Method 1's performance, the 1,000 randomly selected positions $(x_1,y_1,z_1)$ on halo orbits are inputted into the method. Position errors are determined between the true positions $(x_1,y_1,z_1)$ and the positions evaluated via the $(t,A_z)$ pair the method outputs $(x_2,y_2,z_2)$. Each position error is plotted against the values of $(t, A_z)$ determined by Method 1. Of 1,000 points, Method 1 produced solutions for all but $85$.

\begin{figure}[h!]
	\centering\includegraphics[width=0.8\textwidth]{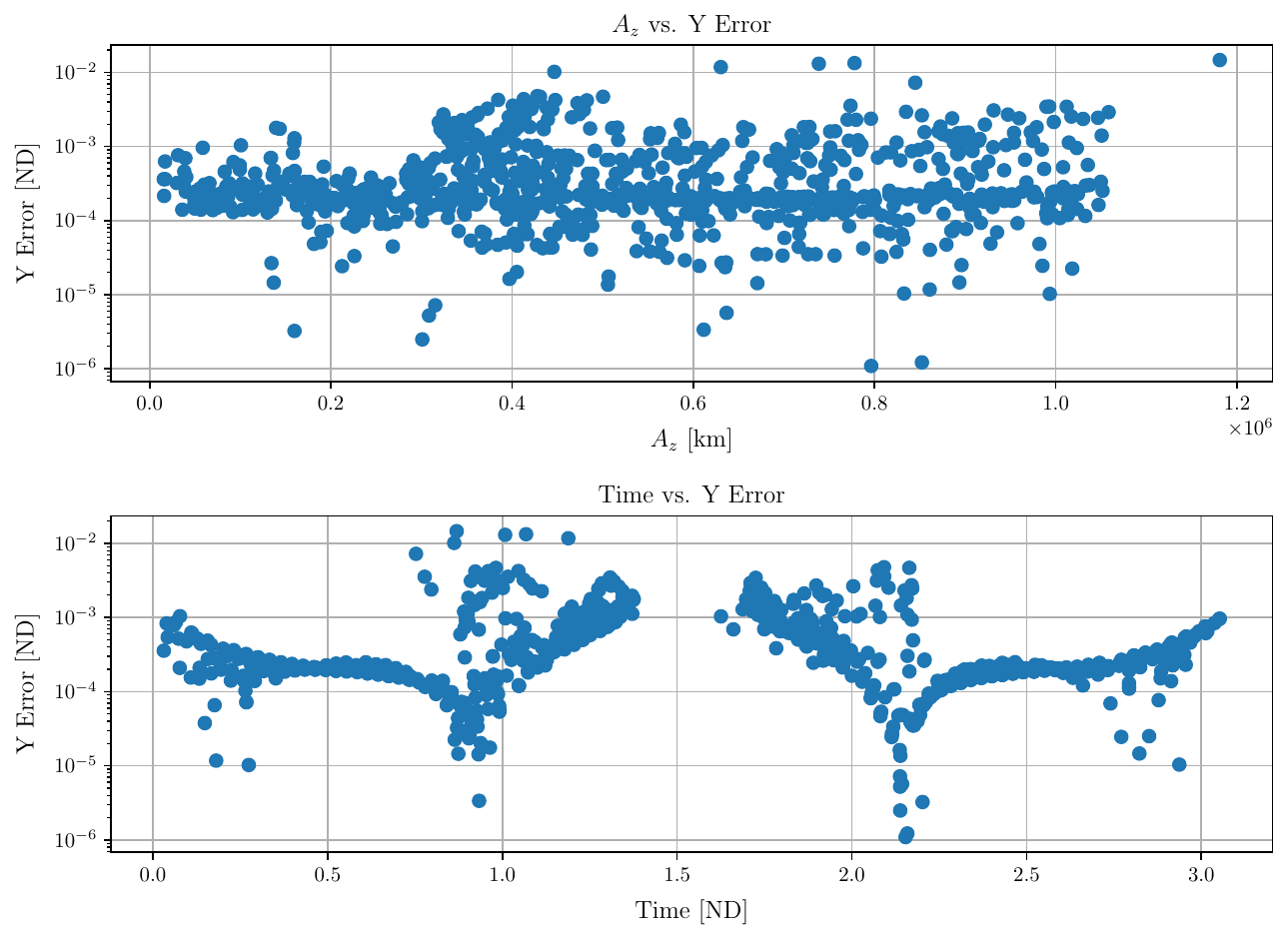}
	\caption{$y$ error compared to time and $A_z$ from Method 1}
	\label{fig:yNorm}
\end{figure}

First, consider the plots for $x$, shown in Figure \ref{fig:xNorm}. There is some indication that higher values of $A_z$ correspond to a larger error in $x$, which suggests that higher $A_z$ values are outside the domain of practical convergence \cite{Jorba1998}. The time plot indicates that errors are higher for points in the second half of the orbit. This is expected due to the \emph{additional} approximate nature of points in the second half of a period discussed previously. The gap around $t = 1.5$ suggests that the process approximates poorly in the region, yielding results that are farther away from $t = 1.5$. This can be explained by the behavior shown in Figure \ref{fig:threeCoordErrors}, which indicates that the $x-$ position error is relatively large in the neighborhood of $t=1.5$. Since Method 1 minimizes $x_1-g(z_1, A_z)$ to find $A_z$, it is evident that the large $x-$ position errors inherent to the Lindstedt-Poincar\'e solutions affect Method 1's accuracy. However, it is worth noting that the position error norm associated with the Lindstedt-Poincar\'e solutions is relatively small around $t=1.5$ (see Figure \ref{fig:lpSol}), suggesting that using the position error norm as an accuracy metric may yield more accurate results.

\begin{figure}[ht!]
	\centering\includegraphics[width=0.8\textwidth]{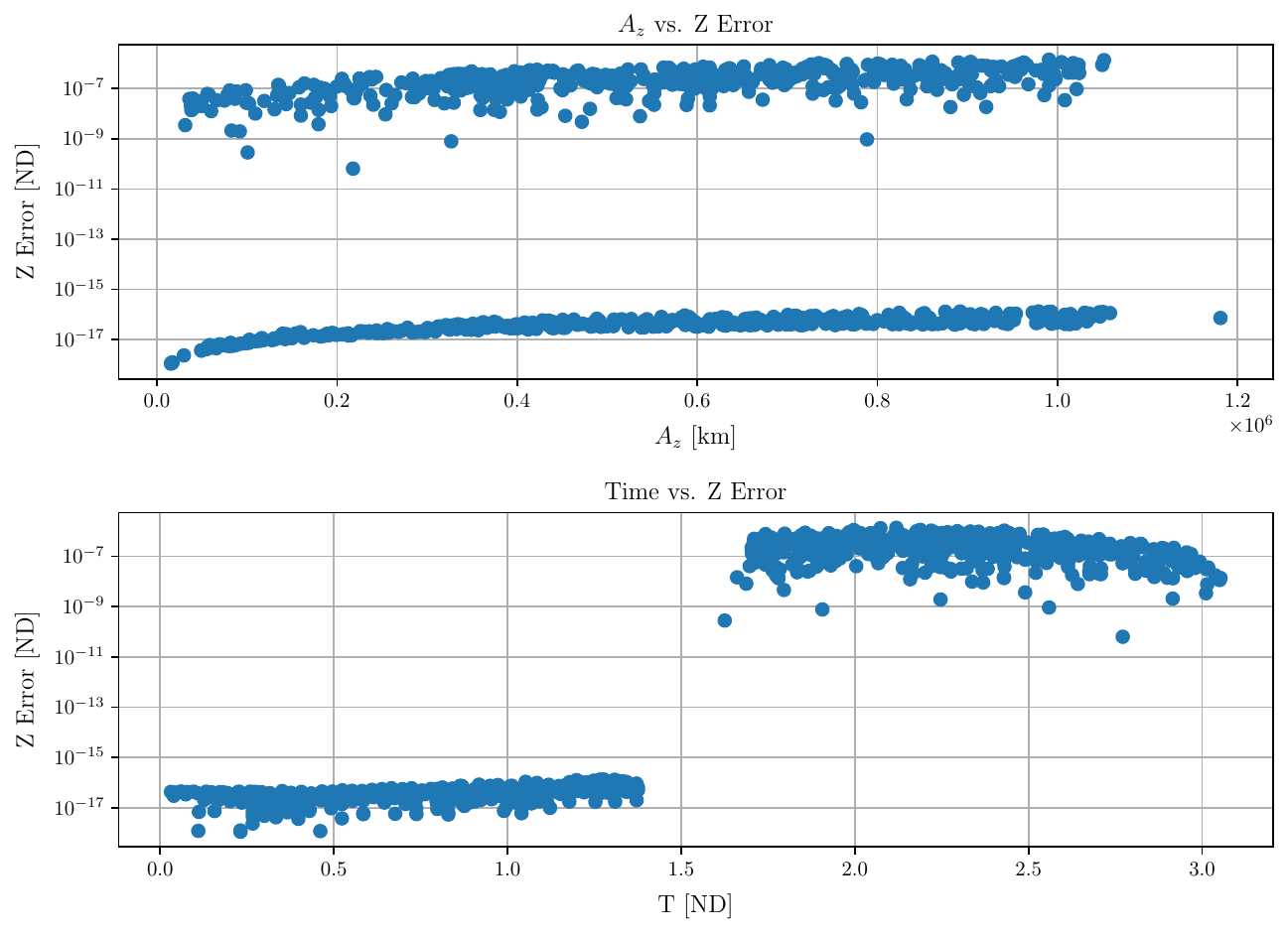}
	\caption{$z$ error compared to time and $A_z$ from Method 1}
	\label{fig:zNorm}
\end{figure}

Figure \ref{fig:yNorm} presents the $y$ error plots. Note that the magnitude of the error is much larger than that of $x$ and $z$. This is expected based on the position error patterns that occur with the Lindstedt-Poincar\'e solutions. Similar to the $A_z$ vs. $x$ error plot, the $y$ error plot suggests that the error increases as $A_z$ increases. It is interesting to note that the $y$ error increases as it approaches the neighborhood of $t=1.5$ from either the left or right. This may be a consequence of the large slope in error around the neighborhood.

The two plots in Figure \ref{fig:zNorm} emphasize two of the conclusions drawn above. The $A_z$ plot clearly shows that the error in $z$ increases as $A_z$ increases. Additionally, the time plot shows that the errors associated with values of time in the second half of the orbit are larger than those of the first half.

\newpage
\subsection{Finding $(t,A_z)$: Method 2}
Method 2 uses Method 1 but performs post-processing that refines Method 1's results. The result $(t,A_z)$ is plugged into the Lindstedt-Poincar\'e solution, yielding $(x_2,y_2,z_2)$. The position error norm is calculated for $(x_2,y_2,z_2)$ with respect to $(x_1,y_1,z_1)$. If the error norm exceeds $10^{-3}$, the  process below is used to find $(t,A_z)$.
\begin{figure}[h!]
	\centering\includegraphics[width=0.8\textwidth]{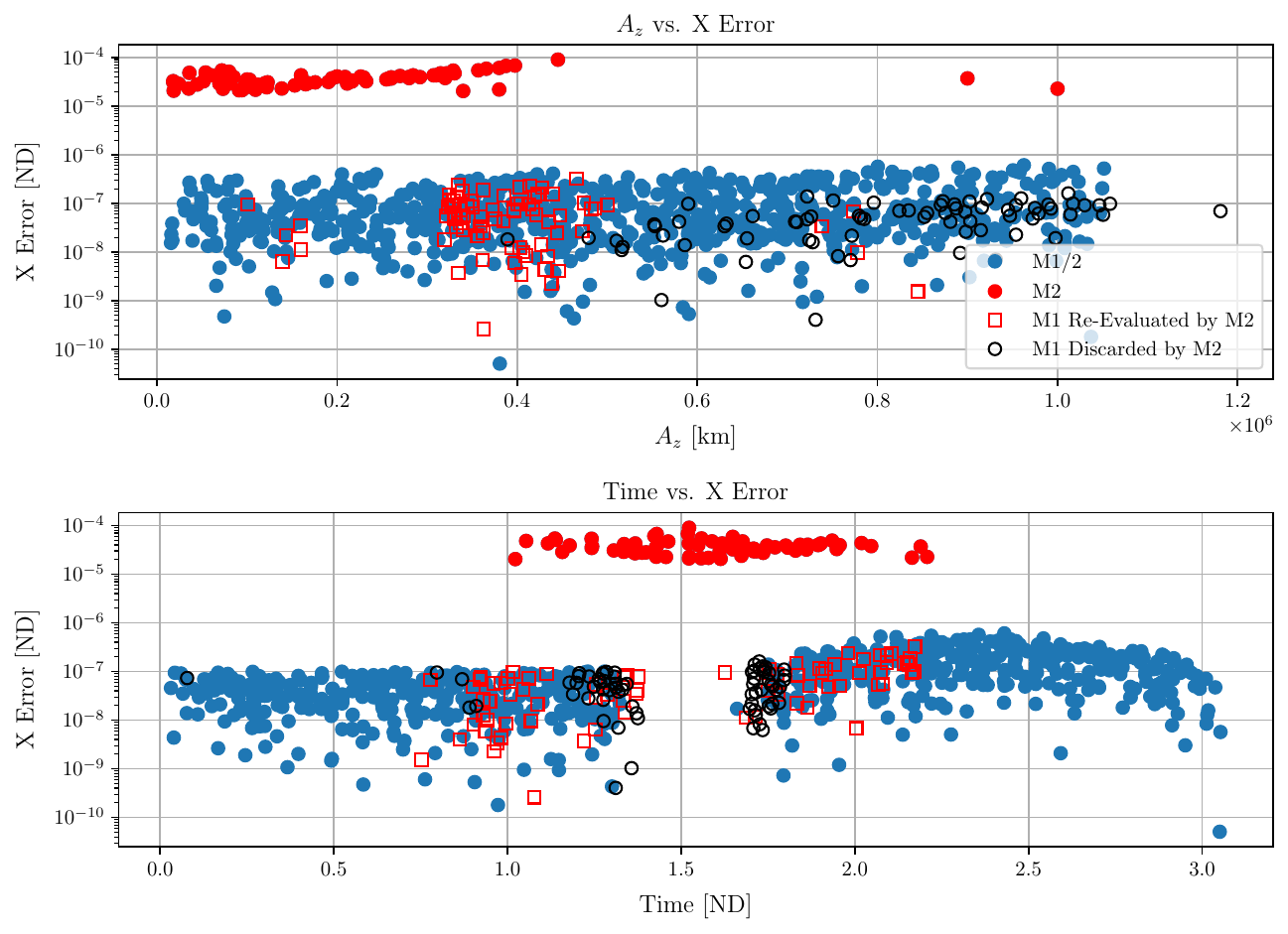}
	\caption{$x$ error compared to time and $A_z$ from Method 2}
	\label{fig:xNorm2}
\end{figure}
\begin{enumerate}
	\item Initialize a range of $A_z$ values between $100$ km and $1.0 \times 10^6$ km.	
	\item  Initialize a variable $n_0$ to a large magnitude.
	\item In a for loop,
	\begin{enumerate}
	\item Plug in $A_z[i]$ in to $t=f(z_1,A_z)$. Using the Lindstedt-Poincar\'e solutions, substitute in $t$ and $A_z$. Find the position norm error of the resulting $(x_2,y_2,z_2)$ compared to $(x_1,y_1,z_1)$. If the resulting position norm error is less than $n_0$, the index $i$ is stored and the value of $n_0$ is updated.
	\end{enumerate}
	\item Upon completing the for loop, the value of $n_0$ is compared against a user-defined tolerance equal to $10^{-4}$ here.
	\begin{enumerate}
	\item If $n_0 \leq tol$, the process has found an acceptable value of $A_z$ corresponding to the stored index. The method then proceeds to find the associated value of $t$.
	\item If $n_0 > tol$, the value of $A_z$ associated with the smallest $n_0$ is used to define a new range of $A_z$ values centered around it. This range gets smaller upon each iteration. If the user-defined number of iterations is exceeded, the process exits and indicates that no $A_z$ value could be found.
	\end{enumerate}
	\item The value of $A_z$ that satisfied the tolerance requirement is plugged into $t = f(z_1,A_z)$. If the sign of $y_1$ is positive, then $t$ is adjusted to be in the second half of orbit as discussed above.
\end{enumerate}

\begin{figure}[h!]
	\centering\includegraphics[width=0.8\textwidth]{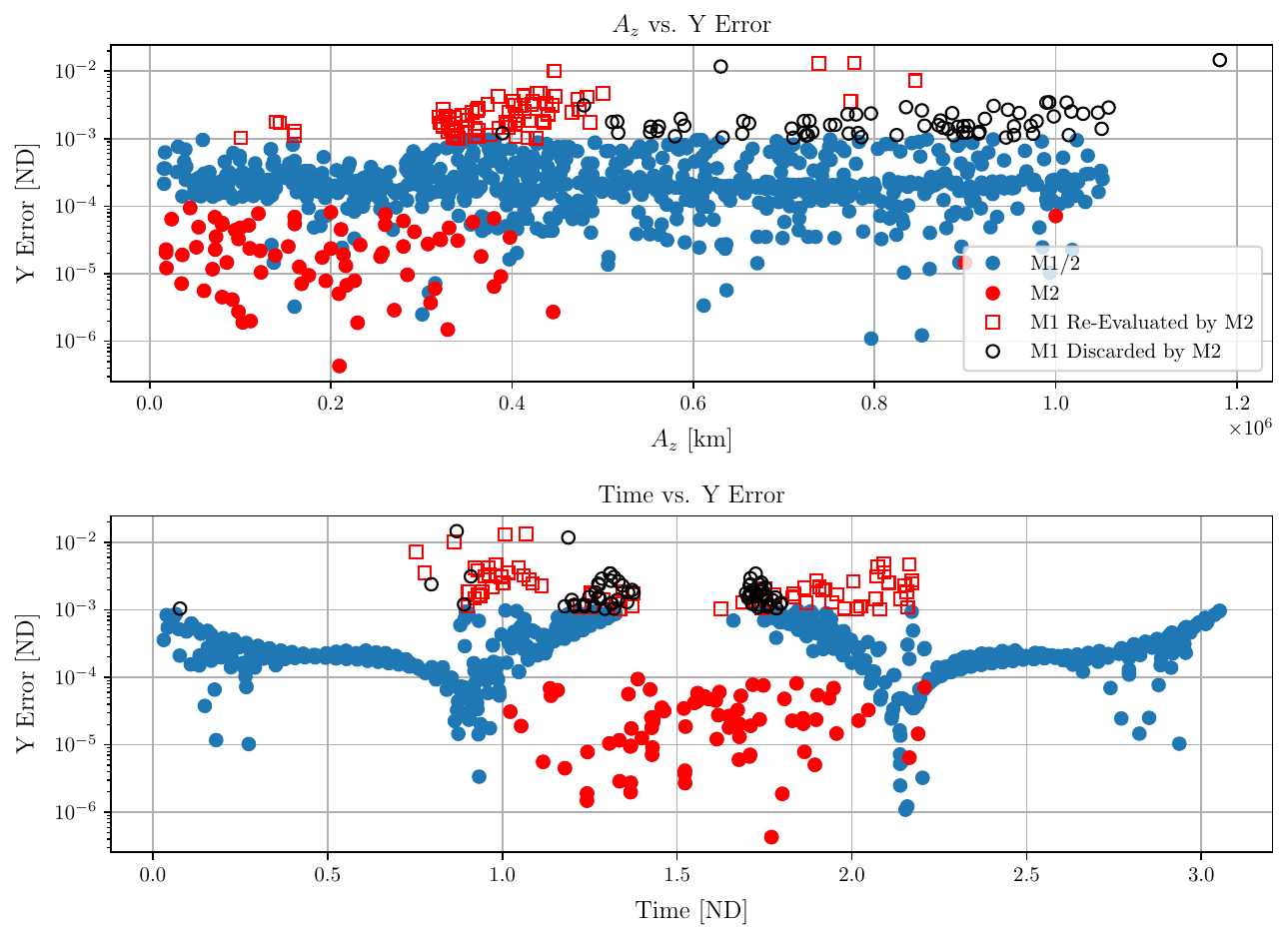}
	\caption{$y$ error compared to time and $A_z$ from Method 2}
	\label{fig:yNorm2}
\end{figure}

The results of Method 2 are either values of $(t, A_z)$ or empty arrays indicating that the method could not find a solution that satisfied the tolerance threshold.

\subsubsection{Results for Method 2}
Since Method 2 refines the solutions $(t, A_z)$ which produced position error norms greater than $10^{-3}$, it is expected that the individual coordinate errors will change and that some solutions will be discarded as unsatisfactory in terms of reaching a desired position error norm less than $10^{-4}$. In the sets of plots for Method 2, additional demarcation of data points is used to distinguish the solutions common to Methods 1 and 2, the solutions unique to Method 2, the solutions from Method 1 that were \textbf{re-evaluated} by Method 2, and the solutions from Method 1 that Method 2 \emph{discarded} as being unsatisfactory. Of 1,000 points, Method 2 produced solutions for all but 156.

Figure \ref{fig:xNorm2} presents the $x$ error for Method 2. Comparing it to Figure \ref{fig:xNorm} indicates that the maximum $x$ error has increased, but there are now solutions in the neighborhood of $t = 1.5$. This shows that a consequence of the solution correction performed by Method 2 is that some of the $x$ errors increased to yield more accurate solutions in the other coordinates. Additionally, it is evident that the majority of solutions discarded from Method 1 are on the edges of the $t=1.5$ neighborhood.

\begin{figure}[h!]
	\centering\includegraphics[width=0.8\textwidth]{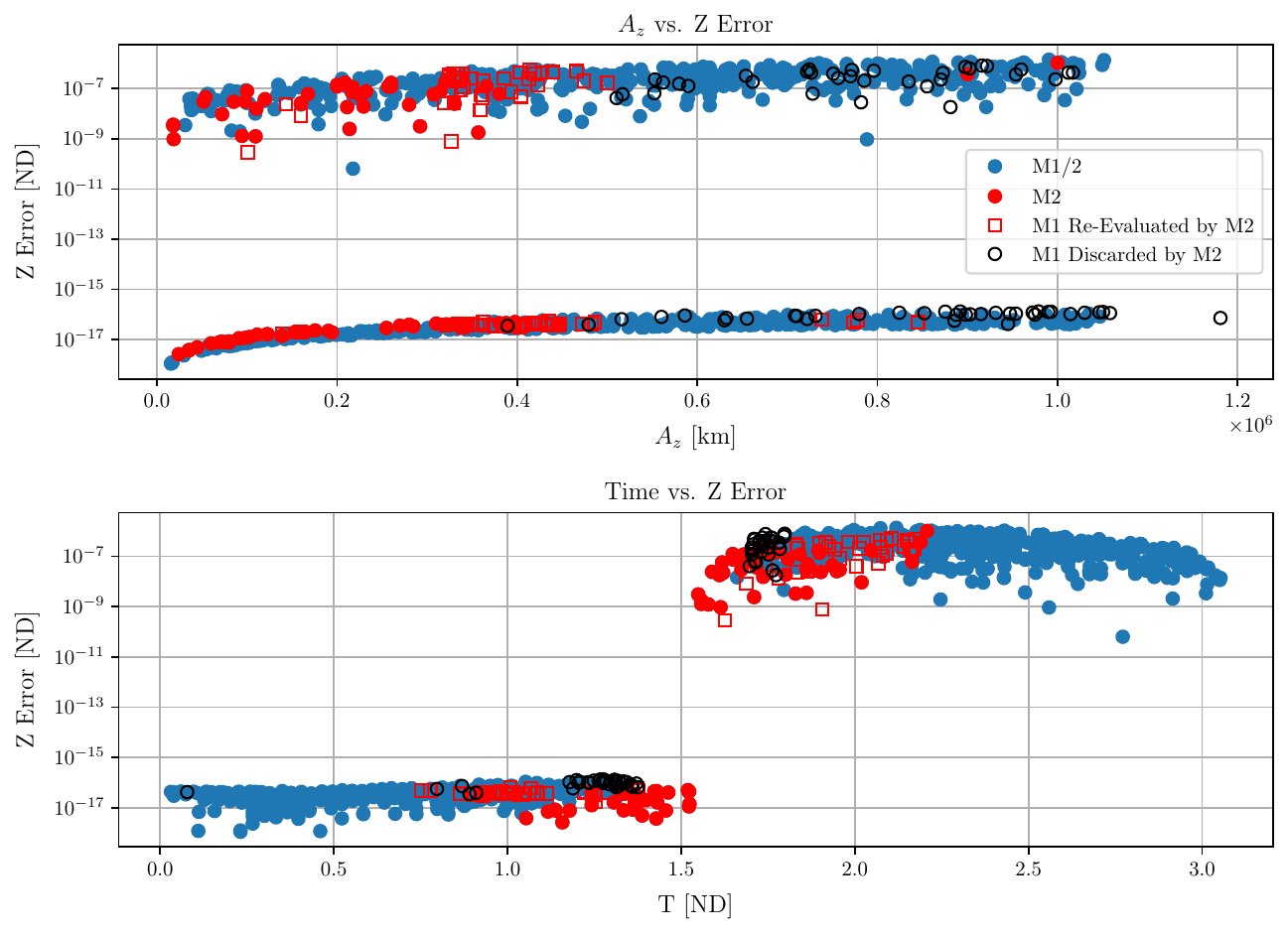}
	\caption{$z$ error compared to time and $A_z$ from Method 2}
	\label{fig:zNorm2}
\end{figure}
Note that the $y$ error, shown in Figure \ref{fig:yNorm2}, has decreased by an order of magnitude using Method 2 compared to Method 1. This is illustrated by the data marked as either re-evaluated or discarded. The time plot indicates that the $y$ errors increase as the solutions approach and depart from $t = 1.5$. Figure \ref{fig:zNorm2} presents the $z$ error using Method 2. The error is comparable to that found using Method 1 (Figure \ref{fig:zNorm}) except that there are now solutions in the neighborhood of $t=1.5$.

\begin{figure}[h!]
	\centering\includegraphics[width=0.8\textwidth]{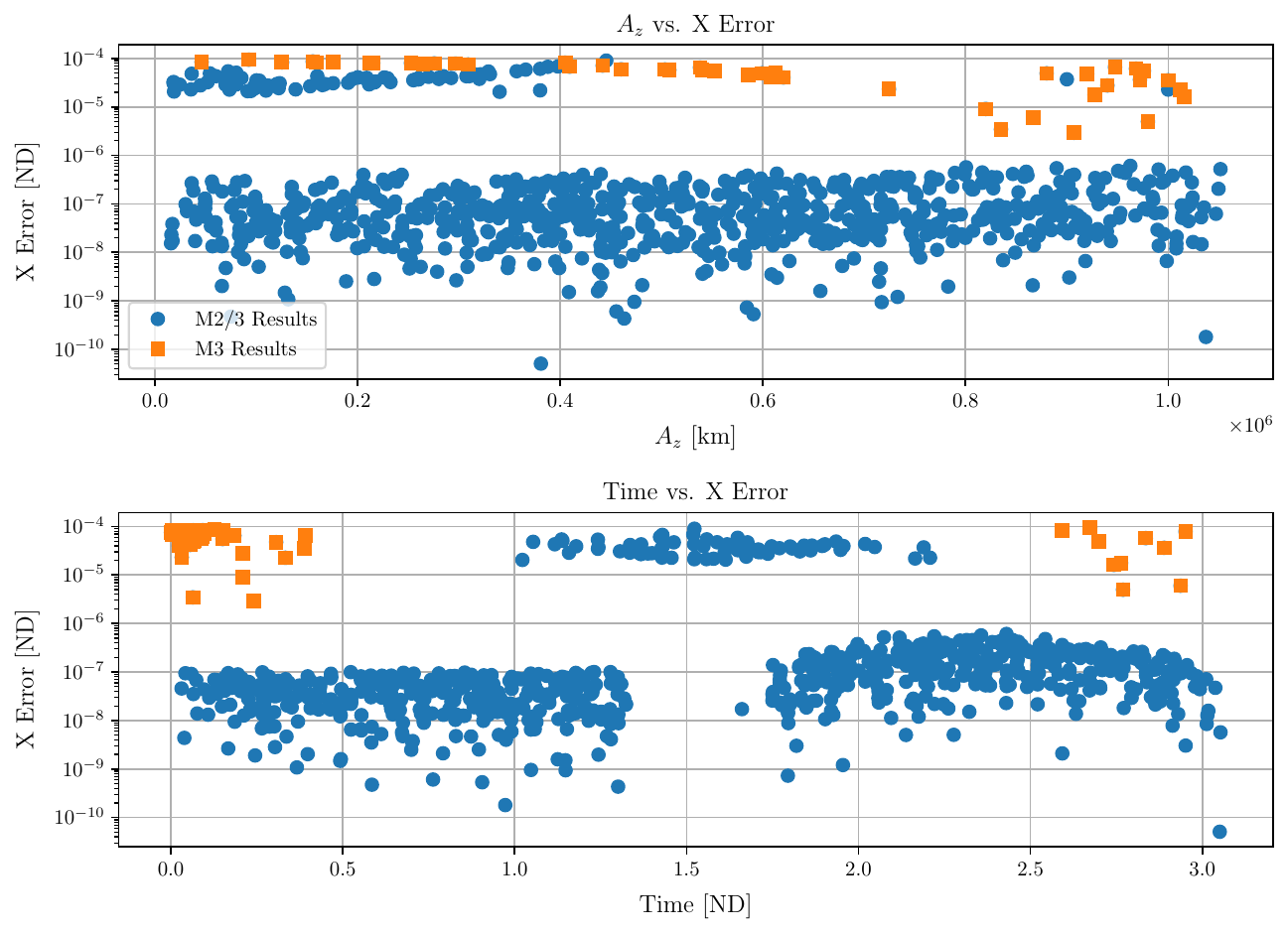}
	\caption{$x$ error compared to time and $A_z$ from Method 3}
	\label{fig:xNorm3}
\end{figure}
\subsection{Finding $(t,A_z)$: Method 3}
Method 3 uses Methods 1 and 2. Specifically, Method 3 executes Method 2; it then considers the points that Method 1 was unable to solve. For those points, it executes steps (1)-(5) described in Method 2. Solutions are considered acceptable if the resulting position error norm is less than $10^{-4}$. The information is saved if solutions for $(t, A_z)$ are found. Otherwise, the method yields empty arrays.

\subsubsection{Results for Method 3}
\begin{figure}[h!]
	\centering\includegraphics[width=0.8\textwidth]{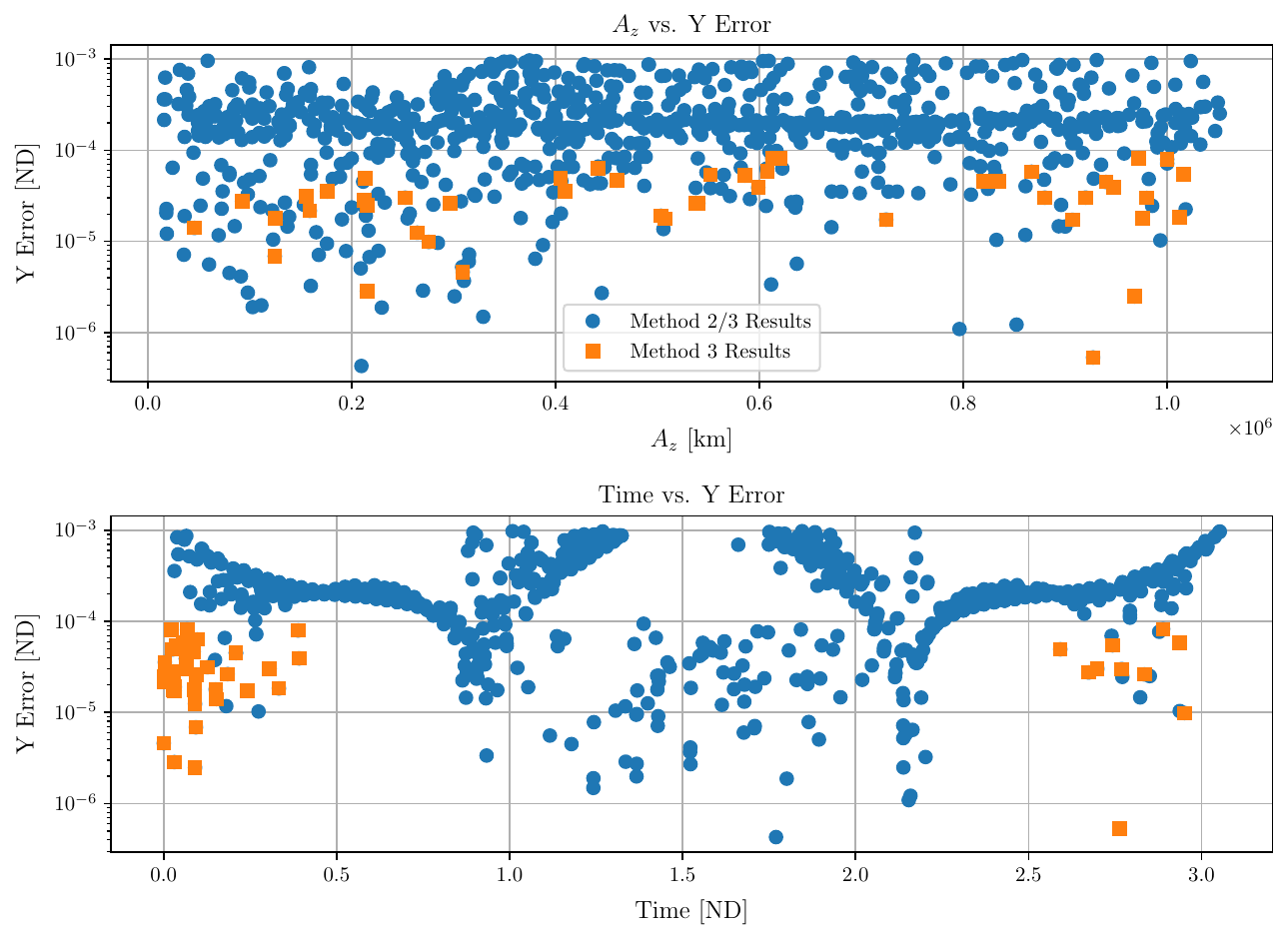}
	\caption{$y$ error compared to time and $A_z$ from Method 3}
	\label{fig:yNorm3}
\end{figure}

\begin{figure}[h!]
	\centering\includegraphics[width=0.8\textwidth]{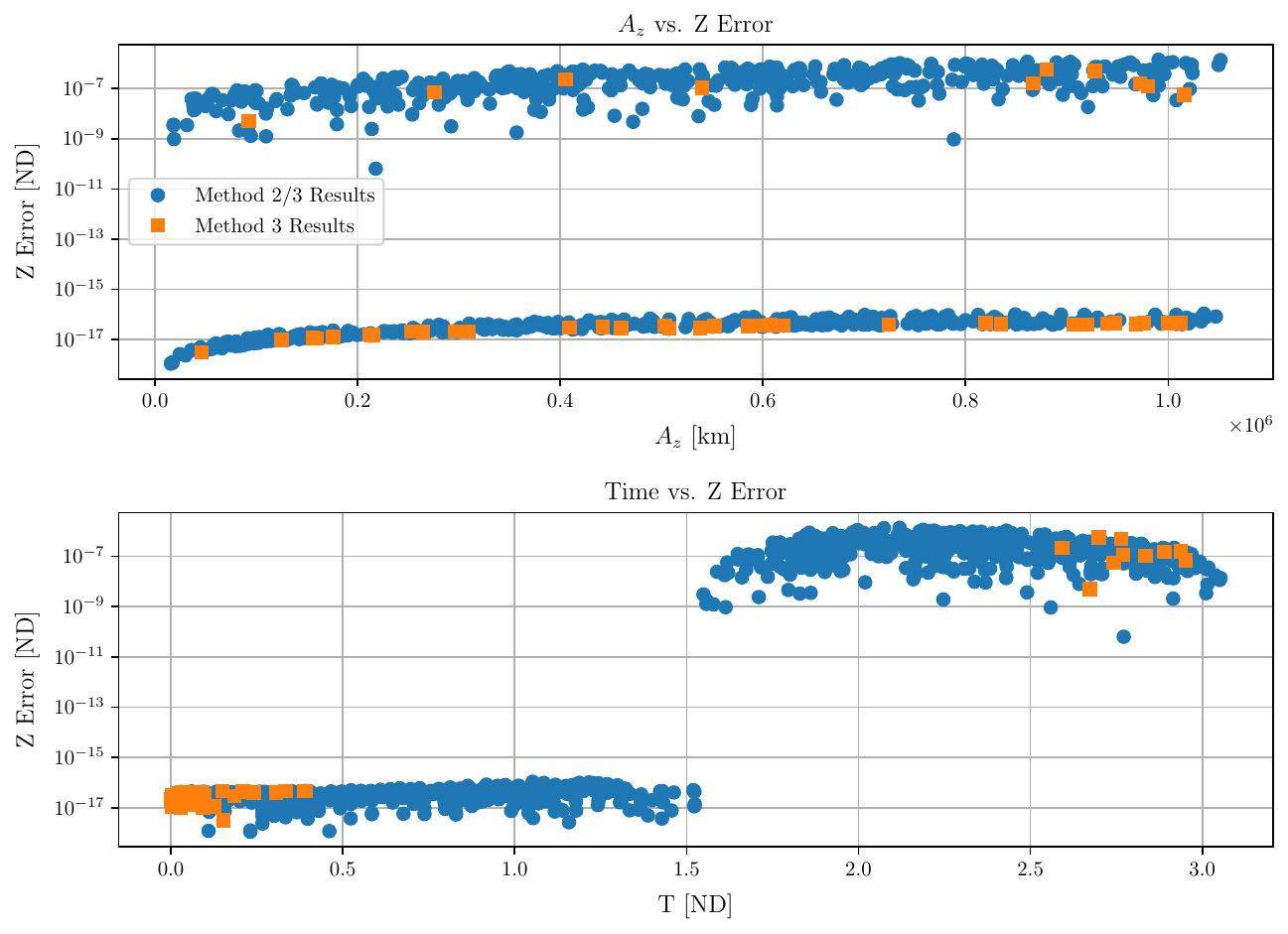}
	\caption{$z$ error compared to time and $A_z$ from Method 3}
	\label{fig:zNorm3}
\end{figure}

Since Method 3 runs Method 2 and then re-evaluates the points that Method 1 was not able to solve within the desired tolerance, the error plots associated with Method 3 are similar to that of Method 2. However, the error plots (Figures \ref{fig:xNorm3}, \ref{fig:yNorm3}, and \ref{fig:zNorm3}) include additional solutions that Method 1 was unable to solve. Of 1,000 points, Method 3 produced solutions for all but 109. Thus, it was able to find 47 more solutions than Method 2. These additional solutions come out of the 85 points that Method 1 could not solve.

Figure \ref{fig:xNorm3} shows the $x$ errors with respect to time and $A_z$. Comparing this figure to the corresponding one for Method 2 (Figure \ref{fig:xNorm2}), it is noted that the new solutions are in the neighborhoods of $t \in [0.0,0.5]$ and $t \in [2.5,3.0]$.
The same is shown in Figures \ref{fig:yNorm3} and Figure \ref{fig:zNorm3}. The location of the additional solutions Method 3 finds shows that the position error norm method sacrifices lower $x$ errors for lower $y$ errors but maintains the same order of magnitude for $z$ errors. 

\section{Conclusion}
The capability to identify the halo orbit that a point resides on in the regime of the restricted three-body problem enables computational efficiency in some mission trajectory designs. Specifically, for a mission like HabEx where there are keepout constraints to consider and several targets to image, this capability enables efficiently finding halo orbits to stitch together rather than employing a blind search using numerical integration of the equations of motion. Using third-order Lindstedt-Poincar\'e solutions, a function of time in terms of the $z$-coordinate and out-of-plane amplitude and a function of the $x$-coordinate in terms of $z$ and the out-of-plane amplitude are found. Three methods that use these functions were developed, each prioritizing a different metric of accuracy. Using 1,000 randomly selected points on known halo orbits, these methods were tested and compared. Ultimately, Method 3 produces the most accurate results, not only in terms of position norm errors but also concerning matching the true time and $A_z$ values and discarding solutions that weighted the performance of the $x$ error at the detriment of the position error norm. Future work includes further investigation into the behavior around $t=1.5$, evaluating an expression for $A_z$ that depends only on $x_1$ and $z_1$ via Mathematica, and determining if fourth-order Lindstedt-Poincar\'e solutions can be used. Whether the process as-is is sufficiently accurate largely depends on the use case. Since it is computationally unlikely to invert Lindstedt-Poincar\'e solutions greater than order $n=4$, use cases should be limited to the domain of practical convergence corresponding to the order of solutions. 

\bibliographystyle{AAS_publication}
\bibliography{biblio.bib}

\begin{thebibliography}{10}

\bibitem{euler1767motu}
L.~Euler, ``De motu rectilineo trium corporum se mutuo attrahentium,''  {\em
  Novi commentarii academiae scientiarum Petropolitanae}, 1767, pp.~144--151.

\bibitem{lagrange1772essai}
J.-L. Lagrange, ``Essai sur le probleme des trois corps,''  {\em Prix de
  l'acad{\'e}mie royale des Sciences de paris}, Vol.~9, 1772, p.~292.

\bibitem{farquhar1968control}
R.~Farquhar, ``The Control and Use of Libration-Point Satellites, Ph. D.
  Dissertation, Dept. of Aeronautics and Astronautics, Stanford University.
  Stanford, CA, 1968,''  1968.

\bibitem{breakwell1979halo}
J.~V. Breakwell and J.~V. Brown, ``The 'halo' family of 3-dimensional periodic
  orbits in the Earth-Moon restricted 3-body problem,''  {\em Celestial
  mechanics}, Vol.~20, No.~4, 1979, pp.~389--404.

\bibitem{richardson1980analytic}
D.~L. Richardson, ``Analytic construction of periodic orbits about the
  collinear points,''  {\em Celestial mechanics}, Vol.~22, No.~3, 1980,
  pp.~241--253.

\bibitem{richardson1980halo}
D.~L. Richardson, ``Halo orbit formulation for the ISEE-3 mission,''  {\em
  Journal of Guidance and Control}, Vol.~3, No.~6, 1980, pp.~543--548.

\bibitem{howell1984three}
K.~C. Howell, ``Three-dimensional periodic halo orbits,''  {\em Celestial
  mechanics}, Vol.~32, 1984, p.~53.

\bibitem{gomez1991study}
G.~G{\'o}mez, A.~Jorba, J.~Masdemont, and C.~Sim{\'o}, ``Study refinement of
  semi-analytical halo orbit theory,''  {\em Final Report, ESOC Contract},
  Vol.~8625, 1991, p.~89.

\bibitem{gomez1991moon}
G.~G{\'o}mez, {\`A}.~Jorba, J.~Masdemont, and C.~Sim{\'o}, ``Moon’s Influence
  on the Transfer from the Earth to a Halo Orbit Around L 1,''  {\em
  Predictability, Stability, and Chaos in N-Body Dynamical Systems},
  pp.~283--290, Springer, 1991.

\bibitem{gomez1993study}
G.~G{\'o}mez, A.~Jorba, J.~Masdemont, and C.~Sim{\'o}, ``Study of the transfer
  from the Earth to a halo orbit around the equilibrium point L 1,''  {\em
  Celestial Mechanics and Dynamical Astronomy}, Vol.~56, No.~4, 1993,
  pp.~541--562.

\bibitem{howell1994transfer}
K.~Howell, D.~Mains, and B.~Barden, ``Transfer trajectories from Earth parking
  orbits to Sun-Earth halo orbits,''  {\em Spaceflight mechanics 1994}, 1994,
  pp.~399--422.

\bibitem{barden1994using}
B.~T. Barden, ``Using stable manifolds to generate transfers in the circular
  restricted problem of three bodies,''  {\em Master degree thesis. West
  Lafayette: School of Aeronautics and Astronautics, Purdue University}, 1994.

\bibitem{farquhar2001flight}
R.~W. Farquhar, ``The flight of ISEE-3/ICE: origins, mission history, and a
  legacy,''  {\em The Journal of the astronautical sciences}, Vol.~49, No.~1,
  2001, pp.~23--73.

\bibitem{domingo1995soho}
V.~Domingo, B.~Fleck, and A.~I. Poland, ``The SOHO mission: an overview,''
  {\em Solar Physics}, Vol.~162, No.~1, 1995, pp.~1--37.

\bibitem{stone1998advanced}
E.~C. Stone, A.~Frandsen, R.~Mewaldt, E.~Christian, D.~Margolies, J.~Ormes, and
  F.~Snow, ``The advanced composition explorer,''  {\em Space Science Reviews},
  Vol.~86, No.~1, 1998, pp.~1--22.

\bibitem{bennett2003microwave}
C.~L. Bennett, M.~Bay, M.~Halpern, G.~Hinshaw, C.~Jackson, N.~Jarosik,
  A.~Kogut, M.~Limon, S.~Meyer, L.~Page, {\em et~al.}, ``The microwave
  anisotropy probe* mission,''  {\em The Astrophysical Journal}, Vol.~583,
  No.~1, 2003, p.~1.

\bibitem{lo2001genesis}
M.~W. Lo, B.~G. Williams, W.~E. Bollman, D.~Han, Y.~Hahn, J.~L. Bell, E.~A.
  Hirst, R.~A. Corwin, P.~E. Hong, K.~C. Howell, {\em et~al.}, ``Genesis
  mission design,''  {\em The Journal of the Astronautical Sciences}, Vol.~49,
  No.~1, 2001, pp.~169--184.

\bibitem{passvogel2010planck}
T.~Passvogel, G.~Crone, O.~Piersanti, B.~Guillaume, J.~Tauber, J.-M. Reix,
  T.~Banos, P.~Rideau, and B.~Collaudin, ``Planck an overview of the
  spacecraft,''  {\em AIP Conference Proceedings}, Vol.~1218, American
  Institute of Physics, 2010, pp.~1494--1501.

\bibitem{pilbratt2010herschel}
G.~Pilbratt, J.~Riedinger, T.~Passvogel, G.~Crone, D.~Doyle, U.~Gageur,
  A.~Heras, C.~Jewell, L.~Metcalfe, S.~Ott, {\em et~al.}, ``Herschel Space
  Observatory-An ESA facility for far-infrared and submillimetre astronomy,''
  {\em Astronomy \& Astrophysics}, Vol.~518, 2010, p.~L1.

\bibitem{wu2012pre}
W.~Wu, Y.~Liu, L.~Liu, J.~Zhou, G.~Tang, and Y.~Chen, ``Pre-LOI trajectory
  maneuvers of the CHANG'E-2 libration point mission,''  {\em Science China
  Information Sciences}, Vol.~55, No.~6, 2012, pp.~1249--1258.

\bibitem{broschart2009preliminary}
S.~B. Broschart, M.-K.~J. Chung, S.~J. Hatch, J.~H. Ma, T.~H. Sweetser, S.~S.
  Weinstein-Weiss, and V.~Angelopoulos, ``Preliminary trajectory design for the
  Artemis lunar mission,''  {\em Advances in the Astronautical Sciences},
  Vol.~135, No.~2, 2009, pp.~1329--1343.

\bibitem{prusti2016gaia}
T.~Prusti, J.~De~Bruijne, A.~G. Brown, A.~Vallenari, C.~Babusiaux,
  C.~Bailer-Jones, U.~Bastian, M.~Biermann, D.~W. Evans, L.~Eyer, {\em et~al.},
  ``The gaia mission,''  {\em Astronomy \& astrophysics}, Vol.~595, 2016,
  p.~A1.

\bibitem{roberts2015early}
C.~Roberts, S.~Case, J.~Reagoso, and C.~Webster, ``Early mission maneuver
  operations for the deep space climate observatory Sun-Earth L1 libration
  point mission,''  {\em AIAA/AAS Astrodynamics Specialist Conference},
  No.~595, 2015, pp.~1--21.

\bibitem{mcnamara2008lisa}
P.~McNamara, S.~Vitale, K.~Danzmann, L.~P. S.~W. Team, {\em et~al.}, ``Lisa
  pathfinder,''  {\em Classical and quantum gravity}, Vol.~25, No.~11, 2008,
  p.~114034.

\bibitem{jia2018scientific}
Y.~Jia, Y.~Zou, J.~Ping, C.~Xue, J.~Yan, and Y.~Ning, ``The scientific
  objectives and payloads of Chang’E- 4 mission,''  {\em Planetary and Space
  Science}, Vol.~162, 2018, pp.~207--215.

\bibitem{gao2019trajectory}
S.~Gao, W.~Zhou, L.~Zhang, W.~LIANG, D.~LIU, and H.~ZHANG, ``Trajectory design
  and flight results for Chang’e 4-relay satellite,''  {\em Scientia Sinica
  Technologica}, Vol.~49, No.~2, 2019, pp.~156--165.

\bibitem{kovalenko2019orbit}
I.~D. Kovalenko and N.~A. Eismont, ``Orbit design for the
  Spectrum-Roentgen-Gamma mission,''  {\em Acta Astronautica}, Vol.~160, 2019,
  pp.~56--61.

\bibitem{gardner2006james}
J.~P. Gardner, J.~C. Mather, M.~Clampin, R.~Doyon, M.~A. Greenhouse, H.~B.
  Hammel, J.~B. Hutchings, P.~Jakobsen, S.~J. Lilly, K.~S. Long, {\em et~al.},
  ``The james webb space telescope,''  {\em Space Science Reviews}, Vol.~123,
  No.~4, 2006, pp.~485--606.

\bibitem{mars_2016}
K.~Mars, ``Gateway,''  Aug 2016.

\bibitem{hall_2020}
L.~Hall, ``What is capstone?,''  Jul 2020.

\bibitem{koon2000dynamical}
W.~S. Koon, M.~W. Lo, J.~E. Marsden, and S.~D. Ross, ``Dynamical systems, the
  three-body problem and space mission design,''  {\em Equadiff 99: (In 2
  Volumes)}, pp.~1167--1181, World Scientific, 2000.

\bibitem{sanchez2020towards}
W.~D. Sanchez, {\em Towards fuel-efficient formation flying of an observatory
  and external occulter at Sun-Earth L2}.
\newblock PhD thesis, Massachusetts Institute of Technology, 2020.

\bibitem{barden1998fundamental}
B.~T. Barden and K.~C. Howell, ``Fundamental motions near collinear libration
  points and their transitions,''  {\em The Journal of the Astronautical
  Sciences}, Vol.~46, No.~4, 1998, pp.~361--378.

\bibitem{soto2020orbital}
G.~J. Soto, {\em Orbital Design Tools and Scheduling Techniques for Optimizing
  Space Science and Exoplanet-Finding Missions}.
\newblock PhD thesis, Cornell University, 2020.

\bibitem{Jorba1998}
A.~Jorba, ``Dynamics in the centre manifold of the collinear points of the
  {Restricted} {Three} {Body} {Problem},''  {\em Physica D: Nonlinear
  Phenomena}, Vol.~132, July 1998, 10.1016/S0167-2789(99)00042-1.

\end{thebibliography}

\end{document}